\documentclass[runningheads]{llncs}
\usepackage{graphicx}
\usepackage{amsmath}
\usepackage{lscape}
\usepackage{makecell}
\def\st{\mathop{\rm s.t.}}

\begin{document}
\title{Integer programming models for the semi-obnoxious p-median problem}
\titlerunning{Models for the semi-obnoxious p-median problem}
%
\author{Nikolaos Ploskas\inst{1}\orcidID{0000-0001-5876-9945} \and
Kostas Stergiou\inst{1}\orcidID{0000-0002-5702-9096}}
\authorrunning{Ploskas \& Stergiou}
%
\institute{University of Western Macedonia, Kozani 50100, Greece\\
\email{nploskas@uowm.gr}, \email{kstergiou@uowm.gr}}
\maketitle              
\begin{abstract}
The p-median problem concerns the location of 
facilities so that the sum of distances between the demand points and their nearest facility is minimized. We study a variant of this classic location problem where minimum distance 
constraints exist both between the facilities and between the facilities and the demand points. This specific type of problem can be used to model situations where the facilities to be located are semi-obnoxious. But despite its relevance to real life scenarios, it has received little attention within the vast literature on location problems. We present twelve ILP 
models for this problem, coupling three formulations of the p-median problem with four formulations of the distance constraints. 
We utilize Gurobi Optimizer v9.0.3 in order to compare these ILP 
models on a large dataset of problems. Experimental results demonstrate that  
the classic p-median model proposed by ReVelle \& Swain 
and the model proposed by Rosing et al. 
are the best performers.

\keywords{Facility location \and p-median problem \and Integer optimization \and Distance constraints.}
\end{abstract}
\section{Introduction}

Facility location problems 
are among the most widely studied problems in OR, AI, computational geometry, and other disciplines. 
Most such 
problems include ``pull'' objectives, where clients wish to have the facility sites (e.g., pharmacies and stores) located close to them. For instance, the {\em p-median problem} concerns the location of {\em p} facilities so that the sum of (weighted) distances between the demand points and their nearest facility is minimized.

Another widely studied class includes problems with ``push'' objectives, where we seek to place the facility sites away from demand points, and in some cases from each other. Such facilities are known as {\em obnoxious} (or undesirable) because they have hazardous properties (i.e., dump sites for chemicals, 
prisons, etc.). In between, we have the class of {\em semi-obnoxious} (or semi-desirable) location problems where the facilities have both desirable and undesirable properties, meaning that both pull and push objectives are present \cite{krarup2002discrete}. In this case, clients wish to have facilities located close to them, but not too close. For example, the residents of a suburb wish to have restaurants and bars situated near them, but not too close because of the noise and traffic associated with such facilities.  

One way to model the requirements that arise when trying to locate semi-obnoxious facilities is by placing minimum allowed distance constraints between facilities and/or between facilities and demand points \cite{krarup2002discrete}. For instance, gas stations (which are typical semi-obnoxious facilities) need to be at a minimum distance from one another and also from residencies, for safety reasons.
At the same time, they must not be placed too far away from residencies for practical reasons. Such a scenario can be captured as a p-median problem with distance constraints.

The p-median problem, which is $\mathcal{NP}$-hard on general networks for an arbitrary p, was originally proposed by Hakimi \cite{hakimi1964optimum,hakimi1965optimum}. In 1970 ReVelle \& Swain \cite{revelle1970central} presented the first IP 
for the p-median problem. They utilized a structure proposed by Balinski \cite{balinski1965integer} in a plant location problem. Despite the large number of variables and constraints of this problem, it has many integer-friendly properties.

In 1979 Rosing et al. \cite{rosing1979p} proposed a formulation 
by using fewer Balinski constraints and adding a constraint used by Efroymson \& Ray \cite{efroymson1966branch} in a plant location problem. This model reduces the number of constraints needed while retaining many of the integer-friendly properties of the ReVelle \& Swain formulation. In 2003 Church \cite{church2003cobra} proposed a new formulation of the p-median problem by building upon the Rosing et al. model and reducing variables based upon a property of ``equivalent assignment conditions''. There are also other works that propose reformulations of the classic p-median problem \cite{cornuejols1980canonical,densham1992strategies}. For a thorough literature review on formulations of the p-median problem, see \cite{church2003cobra}.


Location problems with distance constraints have received little attention within the vast literature on location problems in OR and related areas. 
There are some early works that considered maximum distance constraints between the demand nodes and the facility locations \cite{church1977results,chvatal1979greedy,tansel1982duality}. 
Moon and Chaudhry were the first to systematically study location problems with distance constraints \cite{moon1984analysis}. Among the problems studied was the 
Median$/l/l$, in their terminology, where we seek to minimize the sum of the distances between the clients and their closest facility (as in p-median) subject to constraints that impose lower bounds (hence the $l$) in the distance 
between any two facilities and also between any facility and any client. 
This is the problem that we revisit in this paper, simply calling it {\em p-median with distance constraints}. 

This specific problem has not received attention since then, but other location problems with distance constraints have been studied \cite{chaudhry1986locating,moon1991minimax,berman2008minimum,drezner2018weber,drezner2019planar}. 
In a work that is relevant to ours, Berman \& Huang \cite{berman2008minimum} presented various formulations for the minimum weighted covering location problem 
aiming to minimize the total demand covered subject to the condition that no two facilities are allowed to be closer than a pre-specified distance. 
Interestingly, they were the first to experimentally compare four, increasingly more elaborate, models for the distance constraints in a location problem. 

In this work we make a detailed study of the p-median problem with distance constraints for the first time. Given that the ILP modeling of this problem is based on two decisions: how to model the p-median objective and how to model the distance constraints, we develop and compare 12 different formulations by combining the three main formulations of the p-median problem with the four formulations of the distance constraints presented by Berman \& Huang \cite{berman2008minimum}, extending them to include minimum distance constraints between facilities and clients.

To compare these ILP 
models, we utilize Gurobi Optimizer v9.0.3 and apply the solver on a set of instances generated using as basis the standard p-median benchmark dataset \cite{beasley1985note}. The experimental evaluation 
shows that the classic p-median model proposed by ReVelle \& Swain and the model proposed by Church are the best performers. Regarding the distance constraint formulations, there is no clear winner. An interesting conjencture that can be made based on our results is that simple models that were considered inefficient in the past are quite competitive on our studied problem due to the advancements in ILP solver technology.



\section{Problem definition}
\label{sec:formulation}

We assume that {\em p} homogeneous facilities are to be placed on a weighted network $G=\{V,E\}$ embedded on a plane, where a set {\em I} of demand points (clients) are already located. We also assume that the set {\em J} of nodes in $G$ where facilities can potentially be located is known. Hence, we deal with a discrete/network location problem. The weight $w_{ij}$ on an edge $(i,j) \in E$ denotes the symmetric service cost (typically the distance) between nodes $i$ and $j$. 

Between each pair of facilities $i$ and $j$ there is a distance constraint $Euc_d(i,j) > d_1$ specifying that the Euclidean distance between the nodes where the facilities $i$ and $j$ are located must be greater than $d_1$, where $d_1$ is a real constant. Also, between each facility $i$ and client $k$ there is a distance constraint $Euc_d(i,k) > d_2$ specifying that the Euclidean distance between the node where facility $i$ is located and node where client $k$ is located must be greater than $d_2$, where $d_2$ is a real constant. The constants $d_1$ and $d_2$ are common for all facilities and clients since we consider the homogeneous case in this work.

At this point let us clarify that the Euclidean distance $Euc_d(i,j)$ between two nodes $i$ and $j$ should not be confused with the minimum distance between the two nodes in the network, which is given by the sum of weights of the edges that belong to the shortest path between $i$ and $j$. It is natural to model the costs of service between $i$ and $j$ through the cost of the paths between the two nodes, but we believe that any distance constraint established due to, say, safety reasons, should consider a metric like the Euclidean distance instead. 

The goal is to minimize the sum of distances between each demand point and its nearest facility subject to the satisfaction of all the distance constraints. The existence of distance constraints between clients and facility locations means that, in contrast to the standard p-median problem where all nodes in the network are clients, in our case the set of potential facility sites is discrete from the set of client nodes.

\section{ILP models}
\label{sec:ILPmodels}

In this section, we present various binary ILP 
models for solving the homogeneous p-median problem with distance constraints so that the sum of distances between the demand points and their nearest facility is minimized. 
The following notation is used:
\begin{itemize} 	
  \item $I$: set of demand nodes.
  \item $J$: set of candidate facility sites.
  \item $p$: the number of facilities to be located.
  \item $s_{ij}$: the shortest distance between any two nodes (demand nodes or facility sites).
  \item $d_1$: the minimum distance between each pair of facilities $\left(i, j\right)$, where $i \in J$ and $j \in J$.
  \item $d_2$: the minimum distance between each facility $i \in J$ and all demand nodes.
  \item $x_{j} = 1$ if a facility is located at facility site $j \in J$ and 0 otherwise.
  \item $y_{ij} = 1$ if a demand node $i \in I$ is assigned to facility site $j \in J$ and 0 otherwise.
\end{itemize}

In subsection~\ref{sec:pmedmodels} we present the three p-median formulations used in this work, while subsection~\ref{sec:distancemodels} includes the four distance constraint models.

\subsection{p-median models}
\label{sec:pmedmodels}

We can formulate the p-median model proposed by ReVelle \& Swain \cite{revelle1970central} as:

\begin{eqnarray}
\min & \sum\limits_{i \in I} \sum\limits_{j \in J} s_{ij} y_{ij} & \label{obj1}\\
\st & \sum\limits_{j \in J} x_j = p & \label{con1a}\\
    & \sum\limits_{j \in J} y_{ij} = 1 & \forall i \in I \label{con1b}\\
    & y_{ij} \leq x_j & \forall i \in I, \forall j \in J \label{con1c}\\
    & x_{j} \in \left\{ 0, 1 \right\} & \forall j \in J \label{con1d}\\
    & y_{ij} \in \left\{ 0, 1 \right\} & \forall i \in I, \forall j \in J \label{con1e} 
\end{eqnarray}

The objective function~\ref{obj1} aims at minimizing the sum of the shortest distances between the demand points and their nearest open facility. Constraint~\ref{con1a} specifies that $p$ facilities are to be opened, while Constraint~\ref{con1b} ensures that each demand node will be served by one facility site. Constraint~\ref{con1c}, i.e., the Balinski constraints, guarantees that each demand node will be served by a facility site that is opened. 

This model has $\left | J \right | + \left | I \right | \times \left | J \right |$ variables and $\left | I \right | + \left | I \right | \times \left | J \right | + 1$ constraints. Although it was almost impossible to solve problems of large size using this model in the past, ILP solvers have become significantly faster and capable of solving very large problems in the last two decades.

Rosing et al. \cite{rosing1979p} noticed that the Balinski constraints make the constraint set too large and they could use the Efroymson \& Ray constraints instead:

\begin{equation}
\label{royefr}
\sum\limits_{i \in I} y_{ij} \leq n x_j, \forall j \in J
\end{equation}

This constraint form adds only $\left | I \right |$ constraints to the model, whereas the Balinski form adds $\left | I \right | \times \left | J \right |$ constraints. However, the Efroymson \& Ray constraints are not integer-friendly, requiring many time consuming branch and bound operations. Rosing et al. formulated a model which balances between the use of Balinski constraints and Efroymson \& Ray constraints with the aim to generate a smaller integer-friendly model. The proposed model is the following:

\begin{eqnarray}
\min & \sum\limits_{i \in I} \sum\limits_{j \in J} s_{ij} y_{ij} & \label{obj2}\\
\st & \sum\limits_{j \in J} x_j = p & \label{con2a}\\
    & \sum\limits_{j \in F_i} y_{ij} = 1 & \forall i \in I \label{con2b}\\
    & y_{ij} \leq x_j & \forall i \in I, \forall j \in K_{ir}, \text{ where } i \neq j \label{con2c}\\
    & \sum\limits_{i \in I} y_{ij} \leq n x_j & \forall j \in J \label{con2d}\\
    & x_{j} \in \left\{ 0, 1 \right\} & \forall j \in J \label{con2e}\\
    & y_{ij} \in \left\{ 0, 1 \right\} & \forall i \in I, \forall j \in J \label{con2f} 
\end{eqnarray}
where $r$ is the Balinski constraint rank-cutoff level, $F_i$ is the set of all facility sites except the $p - 1$ furthest sites from client $i$, and $K_{ir}$ is the set of the r-closest sites to client $i$. 

This model is the same as the original formulation except that: (i) Efroymson \& Ray constraints \ref{con2d} have been included for each facility site, (ii) Balinski constraints \ref{con2c} are included for only the r-closest facility sites to each client, and (iii) $p - 1$ assignments to each client are eliminated. This model has $\left | J \right | + \left | I \right | \times \left ( \left | J \right | - p + 1 \right )$ variables and $2 \left | I \right | + r \left | I \right | + 1$ constraints.

Church \cite{church2003cobra} built upon Rosing et al. model and introduced a new formulation by reducing variables based upon a property of ``equivalent assignment conditions''. The proposed model is the following:

\begin{eqnarray}
\min & \sum\limits_{\left (i, j \right) \in A} C_{ij} y_{ij} & \label{obj3}\\
\st & \sum\limits_{j \in J} x_j = p & \label{con3a}\\
    & \sum\limits_{k \leq \left | J \right | - p + 1} \sum\limits_{ \left ( l, j \right ) \in M_{ik}} y_{lj} = 1 & \forall i \in I \label{con3b}\\
    & \sum\limits_{ i \in A_j} y_{ij} \leq \left ( \left | J \right | - p + 1 \right ) x_j & \forall j \in J \label{con3c}\\
    & y_{ij} \leq x_j & \forall i \text{ and } j \in \left\{ O_{ir} \right\}, \text{ where } i \neq j \label{con3d}\\
    & x_{j} \in \left\{ 0, 1 \right\} & \forall j \in J \label{con3e}\\
    & y_{ij} \in \left\{ 0, 1 \right\} & \forall i \in I, \forall j \in J \label{con3f} 
\end{eqnarray}
where:
\begin{itemize}
    \item $k$: index that refers to the closest facility site order.
    \item $O_{ik}$: index of the $k \in J$ closest facility site to a client $i \in I$.
    \item $\left\{ O_{ik} \right\} = \left\{ O_{i1}, O_{i2}, \cdots, O_{ik} \right\}$: the set of facility sites that comprise the $k \in J$ closest facility sites to client $i \in I$.
    \item $a_{iklj} = \left\{\begin{array}{ll}
    1, & \text{if variable } x_{lj} \text{ represents the } k \in J \text{ closest assignment variable } \\
    & \text{for client } i \in I \\ 
    0, & \text{ otherwise}. 
    \end{array}\right.$
    \item $M_{ik} = \left\{ \left ( l, j \right ) \mid x_{lj} \text{ represents the } k \in J \text{ closest facility assignment for client } i \in I\right\}$.
    \item $C_{lj} =  \sum\limits_{i \in I}  \sum\limits_{k \in J} d_{ik} a_{iklj}$.
    \item $A = \left\{ \left ( i, j \right ) \mid \text{ variable } x_{ij} \text{ is included in the model} \right\}$, after variable substitution and eliminating any assignment to the furthest $p - 1$ facility sites for a given client $i \in I$.
    \item $A_j = \left\{ i \mid \left ( i, j \right ) \in A \right\}$.
\end{itemize}

Each constraint in this model mirrors the associated constraint in the Rosing et al. model. The additional notation used to formulate this model is associated with the concept of equivalent variables and their substitution. For every equivalent assignment condition present in the model, a variable is eliminated and a constraint is eliminated when that variable was present in a Balinski constraint. This model is no larger than the Rosing et al. model.

\subsection{Distance constraints}
\label{sec:distancemodels}

We now present four formulations for the minimum distance constraints between each pair facilities as used by Berman \& Huang \cite{berman2008minimum} for the minimum weighted covering location problem. We modify these formulations and use them to model the p-median problem with distance constraints.
Let us introduce two new  sets:
\begin{itemize}
    \item $M = \left\{(i,j) \mid i,j \in J, Euc_d(i,j) \leq d_1\right\}, \forall i \in J, \forall j \in J$: the set of pairs $(i,j)$ s.t. facilities cannot be in facility sites $i$ and $j$ because $i$ and $j$ are not in a safe distance between each other.
    \item $N = \left\{(j) \mid j \in J, \exists k \in I, Euc_d(j,k) \leq d_2\right\}, \forall j \in J$: a facility cannot be placed in facility site $j$ because there exists a demand node $k$ that is not in safe distance from $j$.
\end{itemize}

Using this notation the first, and simplest, distance constraint that can be added to the p-median models is the following:

\begin{equation}
\label{dist1}
x_i + x_j \leq 1, \forall \left ( i, j \right ) \in M
\end{equation}

Constraint~\ref{dist1} states that two facilities cannot be established if the distance between them is less than $d_1$. 

Another possible formulation for the distance constraint between two facilities is the following:

\begin{equation}
\label{dist2}
p \left ( 1 - x_j \right ) \geq \sum\limits_{k \in Q_j} x_k, \forall j \in M
\end{equation}
where $Q_j = \left\{ k \in M \mid Euc_d(k,j) \leq d_1, k \neq j, \forall j \in M \right\}$. 

If a facility site is opened at $j \in M$, $x_k$ must take a value of $0$ for all $k \in Q_j$. When $x_j = 0$, at most $p$ facilities can be sited in $Q_j$. This constraint was originally proposed
by Moon \& Chaudhry~\cite{moon1984analysis} to formulate the anti-covering location problem. Constraint set~\ref{dist2} has fewer constraints than constraint set~\ref{dist1}.

The third formulation for this distance constraint is the following:

\begin{equation}
\label{dist3}
n_j \left ( 1 - x_j \right ) \geq \sum\limits_{k \in Q_j} x_k, \forall j \in M
\end{equation}
where $n_j$ is the minimum coefficient necessary to impose spatial restrictions between members of the set $Q_j$. In order to calculate $n_j$ a complementary problem needs to be solved. Let us denote the graph $G' = \left ( N', L' \right )$, where $N' = M$ and $L' = \left\{ \left ( j, k \right ) \in M \mid Euc_d(j,k) \leq d_1, j \neq k \right\}$. Let $G^{'}_j$ be the subgraph of $G^{'}$ induced by $Q_j \cup {j}$ and $q_j$ is the cardinality of the maximum independent set of $G^{'}_j$. Although the maximum independent set problem is $\mathcal{NP}$-complete, $G^{'}_j$ is usually a small graph. In this work we calculate $q_j$ with the second version of the Bron \& Kerbosch algorithm \cite{bron1973algorithm}. Once $q_j$ is calculated, the value of $n_j$ can be determined by $n_j = \min \left ( p, q_j \right )$. Constraint~\ref{dist3} is tighter than constraint~\ref{dist2}.

The fourth formulation for the distance constraint between two facility sites is the following:

\begin{eqnarray}
\sum\limits_{k \in H_j} x_k \leq 1, \forall H_j \in L_1 \label{dist41}\\
n^{'}_j \left ( 1 - x_j \right ) \geq \sum\limits_{k \in Q^{'}_j} x_k, \forall Q^{'}_j \in L_2 \label{dist42}
\end{eqnarray}
where $H_j$ is a maximal clique which includes node $j$ of $G'$, $Q^{'}_j$ is the subset of $Q_j \cup j$ after removing all elements in $H_j$. $L_1$ is the set of all $H_j \left ( j \in M \right )$ and $L_2$ the set of all nonempty $Q^{'}_j \left ( j \in M \right )$. In order to obtain $n^{'}_j$ the following ILP problem needs to be solved:

\begin{eqnarray}
\min & \sum\limits_{ v \in Q^{'}_j} x_v & \\
\st & x_u + x_v \leq 1 & \forall u, v \in Q^{'}_j \text{ and } \left ( u, v \right ) \in G' \\
    & x_v \in \left\{ 0, 1 \right\} & \forall v \in Q^{'}_j
\end{eqnarray}

Since $Q^{'}_j$ is usually small, the solution of the ILP model needs only a very small amount of time compared to the solution of the p-median problem.

Finally, we also consider a distance constraint between each facility and each client, where all clients should be at a safe distance ($d_2$) for all opened facilities. This constraint can be formulated as:

\begin{equation}
x_j = 0, \forall j \in N
\end{equation}

This constraint set ensures that each demand node is on a safe distance from each facility site if it is opened. Instead of setting variables $x_{j}$, where $\left( j \right) \in N$, equal to $0$, these variables can be removed from the model.

Coupling the three formulations of the p-median problem presented in subsection~\ref{sec:pmedmodels} with the four formulations of the distance constraints presented in this subsection, we get twelve different formulations of the p-median problem with distance constraints.

\section{Computational study}
\label{sec:experiments}

As this is the first time that a computational study of the p-median problem with distance constraints is performed, there is a lack of benchmarks for this problem, and location problems with distance constraints in general. Hence, we experimented with random instances. We first describe the random generation model, then we give results 
from the twelve formulations. Computations were performed on an Intel i7 CPU 8700 with 16 GB of main memory, a clock of 3.2 GHz, an L1 cache of 348 KB, an L2 cache of 2 MB, and an L3 cache of 12 MB, running under CentOS 8.4. The ILP models were solved using Gurobi v9.0.3.

\subsection{Random generation model}

The random generator takes as input the p-median benchmark dataset \cite{beasley1985note}, which consists of $40$ problems with a set $V$ of $100$--$900$ nodes and a value of $p$ between $5$ and $200$. We present results from twelve selected problems that were reasonably hard for the presented formulations. In the classical p-median problem, all nodes are considered to be clients. However, in our case, there is a minimum distance constraint between facilities and clients. Therefore, we select $\left | J \right | < \left | V \right |$ nodes to be candidate facilities, while the remaining nodes $\left | I \right | = \left | V \right | - \left | J \right |$ are clients. 
We have considered two cases:
\begin{enumerate}
    \item We start by randomly selecting $80$\% of the nodes as candidate facility sites and the remaining $20$\% as clients. In case the resulting candidate facility sites is less than or equal to $p$, then we progressively increase the number of candidate facility sites until  $\left | J \right | > p$. Also, in case the generated instances are infeasible, we again increase the number of candidate facilities (while decreasing the number of clients) until we get feasible instances.
    \item We start by selecting $20$\% of the nodes as candidate facility sites and the remaining $80$\% as clients. Similarly to the previous case, we progressively increase the number of candidate facility sites until $\left | J \right | > p$ and the generated instances are feasible.
\end{enumerate}

In order to select the parameter $d_1$, we find the minimum and maximum distance between all pairs of candidate facility sites and we set $d_1$ equal to a randomly generated number in the range $\left [\min, \min + (\max - \min) / 10 \right ]$. This interval generated feasible problems in most cases. We perform a similar procedure to select the parameter $d_2$. In this case, we calculate the minimum and maximum distance between all candidate facility sites and all clients.

Finally, the p-median benchmark problems include only the distance for each edge in the graph. We used Dijkstra's algorithm to calculate the shortest path $s_{ij}, i \in V, j \in V$ between each pair of nodes. We set the Euclidean distance between each pair of nodes to a random number in the interval $\left [ s_{ij} / 2, s_{ij} \right ]$.

\subsection{Experimental results}

Tables~\ref{tab:1} - \ref{tab:2} present the execution times of the twelve ILP formulations for cases 1 and 2 described above. 
We use the following abbreviations: 
(i) RS1 - RS4: the ReVelle \& Swain formulation using the distance constraints in~\ref{dist1}, \ref{dist2}, ~\ref{dist3}, and ~\ref{dist41} - \ref{dist42}, respectively, (ii) RRR1 - RRR4: the Rosing et al. formulation using the 
constraints in~\ref{dist1}, \ref{dist2}, ~\ref{dist3}, and ~\ref{dist41} - \ref{dist42}, respectively, and (iii) CHU1 - CHU4: the Church formulation using the 
constraints in~\ref{dist1}, \ref{dist2}, ~\ref{dist3}, and ~\ref{dist41} - \ref{dist42}, respectively. Results from the best method in each class are given in bold. 

All formulations were executed on ten instances for each of the problems. Tables~\ref{tab:1} - \ref{tab:2} present the average execution times. As the Rosing et al. and Church formulations require as input the parameter $r$, i.e., the Balinski constraint rank-cutoff level, we ran these models for values of $r$ starting from $10$ to $60$ with a step of $10$. As there is no value for $r$ that gives the best results on all problem classes, Tables~\ref{tab:1} - \ref{tab:2} report the results for $r = 50$, which resulted in the best results on most classes. But in the two rows giving the average execution times we also include, in parentheses, the average execution times obtained by running the two models on each problem class with the best $r$ value for the specific class.

The average results show that the ReVelle \& Swain formulation is the best one for problem of Case 1 (Table~\ref{tab:1}). This is contrary to the literature, since computational studies on p-median formulations concluded that the ReVelle \& Swain cannot be used in practice due to its size \cite{church2003cobra}. However, 
current state-of-the-art ILP solvers have become much faster and are able to solve very large instances, using sophisticated preprocessing routines to greatly reduce the size of the problems. The main reason for the superiority of the classic p-median problems in the first case is that it performs much fewer simplex iterations that the other p-median formulations. 

When we customize the value of $r$ for each problem class, the Church formulation is the best one on problems of Case 2 (Table~\ref{tab:2}). However, the other two formulations are not far off. But given that the customization of $r$ 
would require extensive experimentation, the 
results from the case where $r$ is set to a good overall value (50 in this case) give a more accurate picture.

\begin{landscape}
\begin{table}[htbp]
  \centering
  \caption{Comparison of ILP formulations for problems where $\left | J \right | > \left | I \right |$.} \label{tab:1}
    \begin{tabular}{c|cccc|cccc|cccc}
    \hline
    Problem \textbackslash{} Model & RS1   & RS2   & RS3   & RS4   & RRR1  & RRR2  & RRR3  & RRR4  & CHU1  & CHU2  & CHU3  & CHU4 \\
    \hline
    pmed21 & 9.11  & 9.52  & 10.64 & 10.32 & 10.59 & 9.92  & 11.03 & \textbf{8.51} & 9.98  & 9.03  & 11.10 & 8.61 \\
    pmed22 & 9.39  & 8.54  & 10.04 & 7.68  & 7.78  & 7.54  & 8.36  & 8.47  & 7.63  & 7.64  & \textbf{7.34} & 8.03 \\
    pmed23 & 0.90  & 1.89  & 2.17  & 1.04  & 0.61  & 1.00  & 1.07  & 0.67  & 0.68  & 1.17  & 1.03  & \textbf{0.53} \\
    pmed26 & 7.86  & 10.61 & \textbf{7.29} & 8.10  & 9.82  & 11.30 & 12.20 & 9.17  & 10.45 & 11.26 & 8.01  & 12.16 \\
    pmed27 & 6.86  & 7.92  & 6.27  & 6.71  & 5.35  & 8.18  & 6.98  & \textbf{4.52} & 6.29  & 5.72  & 5.95  & 5.60 \\
    pmed28 & 0.55  & 0.86  & 0.81  & 0.78  & \textbf{0.22} & 0.45  & 0.43  & 0.31  & 0.28  & 0.62  & 0.33  & 0.31 \\
    pmed31 & 55.19 & 52.44 & \textbf{50.56} & 53.23 & 60.94 & 63.10 & 64.65 & 80.31 & 60.14 & 56.92 & 55.25 & 54.36 \\
    pmed32 & 9.58  & 12.87 & 11.43 & 11.42 & 12.91 & 12.69 & 12.77 & \textbf{8.51} & 10.80 & 9.40  & 10.34 & 11.68 \\
    pmed33 & 1.25  & 2.20  & 2.88  & 1.95  & \textbf{0.63} & 1.47  & 1.13  & 0.78  & 0.75  & 1.37  & 1.41  & 0.81 \\
    pmed36 & 38.63 & \textbf{37.53} & 39.54 & 40.12 & 48.31 & 38.40 & 40.67 & 43.15 & 42.67 & 37.89 & 38.53 & 40.15 \\
    pmed37 & 2.99  & 5.61  & 5.58  & 3.64  & \textbf{1.85} & 3.76  & 4.63  & 2.85  & 2.21  & 4.44  & 3.62  & 2.70 \\
    pmed38 & \textbf{19.60} & 20.74 & 20.13 & 20.43 & 86.73 & 88.54 & 69.73 & 85.53 & 86.01 & 94.56 & 72.81 & 87.88 \\
    \hline
    Average & \textbf{13.49} & 14.23 & 13.94 & 13.79 & \makecell{20.48\\(18.50)} & \makecell{20.53\\(16.20)} & \makecell{19.47\\(17.27)} & \makecell{21.07\\(17.09)} & \makecell{19.82\\(15.38)} & \makecell{20.00\\(15.54)} & \makecell{17.98\\(15.94)} & \makecell{19.40\\(15.22)} \\
    \hline
    \end{tabular}%

\caption{Comparison of ILP formulations for problems where $\left | I \right | \geq \left | J \right |$.} \label{tab:2}
    \begin{tabular}{c|cccc|cccc|cccc}
    \hline
    Problem \textbackslash{} Model & RS1   & RS2   & RS3   & RS4   & RRR1  & RRR2  & RRR3  & RRR4  & CHU1  & CHU2  & CHU3  & CHU4 \\
    \hline
    pmed21 & 4.56  & 5.93  & 4.58  & 4.29  & 4.36  & 6.25  & 6.24  & 6.20  & 4.75  & \textbf{4.28} & 4.35  & 5.88 \\
    pmed22 & 7.43  & 7.80  & 6.27  & 9.34  & 7.57  & 8.01  & 6.59  & 10.52 & 8.83  & 7.31  & \textbf{5.84} & 7.79 \\
    pmed23 & 0.38  & 0.44  & 0.38  & 0.36  & 0.38  & 0.36  & 0.36  & 0.39  & 0.38  & 0.34  & \textbf{0.31} & 0.32 \\
    pmed26 & 25.15 & 23.92 & 26.58 & \textbf{20.49} & 24.03 & 21.56 & 23.41 & 30.76 & 21.26 & 27.61 & 24.99 & 21.61 \\
    pmed27 & 17.74 & 20.47 & 23.43 & \textbf{17.33} & 25.29 & 21.49 & 24.59 & 24.48 & 18.36 & 24.41 & 24.44 & 28.77 \\
    pmed28 & 0.67  & 0.82  & 0.88  & 0.70  & \textbf{0.52} & 0.78  & 0.68  & 0.62  & 0.62  & 0.68  & 0.54  & 0.57 \\
    pmed31 & 55.19 & 52.44 & \textbf{50.56} & 53.23 & 60.94 & 63.10 & 64.65 & 80.31 & 60.14 & 56.92 & 55.25 & 54.36 \\
    pmed32 & 88.10 & 72.91 & 86.84 & 91.02 & 57.98 & 50.66 & 57.51 & 55.55 & 49.06 & 70.38 & 55.19 & \textbf{44.90} \\
    pmed33 & 1.00  & 1.19  & 1.13  & 0.99  & 0.75  & 1.09  & 1.14  & 0.73  & \textbf{0.68} & 1.00  & 0.89  & 0.73 \\
    pmed36 & 51.52 & 60.16 & 42.44 & 54.91 & 28.63 & 30.34 & 26.89 & 25.44 & 31.45 & \textbf{24.32} & 24.73 & 30.26 \\
    pmed37 & 1.69  & 1.91  & 2.34  & 1.63  & 1.53  & 1.95  & 1.62  & 1.21  & \textbf{1.21} & 1.39  & 1.38  & 1.26 \\
    pmed38 & 146.40 & 144.00 & \textbf{111.81} & 180.86 & 149.10 & 116.81 & 165.18 & 168.31 & 170.86 & 172.81 & 139.64 & 168.68 \\
    \hline
    Average & 33.32 & 32.67 & 29.77 & 36.26 & \makecell{30.09\\(24.72)} & \makecell{\textbf{26.87}\\(23.41)} & \makecell{31.57\\(25.94)} & \makecell{33.71\\(24.95)} & \makecell{30.63\\(23.32)} & \makecell{32.62\\(24.63)} & \makecell{28.13\\(23.13)} & \makecell{30.43\\(26.89)} \\
    \hline
    \end{tabular}%
\end{table}%

\end{landscape}

Case 2 problem instances 
are harder on average and the Rosing et al. formulation is the best performer in this case. Case 2 instances is harder than Case 1 because they are much larger in size (both in terms of variables and constraints), thus leading to more explored nodes in the branch-and-bound tree and more simplex iterations.

Overall, if we consider only relatively hard problem instances that all formulations need more than $10$ seconds to solve, then RS3 is $1.51$ - $1.90$ times faster than all Rosing et al. and Church formulations in Case 1 problems, whereas RRR2 is $1.07$ - $1.37$ times faster than all ReVelle \& Swain and Church formulations in Case 2 ones. Note that some of the generated instances are quite hard. The most difficult instance of Case 1 (resp. Case 2) required $1,548$ (resp. $1,666$) seconds.

Regarding the distance constraints, there is no clear winner, indicating that the more complex formulations do not always pay off. Most best execution times are achieved by the distance constraints in~\ref{dist1},~\ref{dist3} and~\ref{dist41}-~\ref{dist42}. The distance constraint in~\ref{dist3} leads to better average execution times with the Rosing et al. and Church formulations both in Case 1 problems and with the ReVelle \& Swain and Church formulations in Case 2 problems. Therefore, distance constraint in~\ref{dist3} can be regarded as the best candidate for most cases on average.

\vspace{-2mm}

\section{Conclusion}
\label{sec:conclusion}

\vspace{-2mm}

We have investigated for the first time a variant of the p-median problem which concerns the location of semi-obnoxious facilities in the presence of distance constraints between facilities and also between facilities and demand points. 
We have presented twelve ILP models, coupling three formulations of the p-median problem with four formulations of the distance constraints. Experiments have demonstrated that, perhaps surprisingly, the classic p-median formulation coupled with a simple model for the distance constraints is quite competitive as it is the best one in problems with more candidate facility sites than clients, and it is not far off the the best performer (i.e. the Rosing et al.) in problems with more clients than candidate facility sites.

As future work we would like to consider specialized branch-and-bound algorithms for this problem and extend our study 
to the case where the facilities are not homogeneous. 
We can also extend/modify the proposed techniques to related types of location problems, and of course we can investigate the applicability of the methods to real-world problems.

%
%
%
\bibliographystyle{splncs04}
\bibliography{semihomogeneousbib}

\end{document}